# APPROXIMATE SOLUTION OF INVERSE BOUNDARY PROBLEM FOR THE HEAT EXCHANGE PROBLEM BY A.N. TIKHONOV'S REGULARIZATION METHOD

## V.F. Mirasov[1], A.I. Sidikova[2]

[1] **Mirasov Vadim Faritovich** - Post-graduate student, Calculating Mathematics Department, South Ural State University, Chelyabinsk, Russia

[2] **Sidikova Anna Ivanovna** – PhD. (Physics and Mathematics), associate professor, Calculating Mathematics Department, South Ural State University, Chelyabinsk, Russia

Corresponding author: vadim.mirasov@yahoo.com

In this paper the approximate solution of the heat exchange problem by A.N. Tikhonov's regularization method is presented. The errror estimation of this approximate is obtained.

**The direct problem definition**

Let the thermal process is described by the equation

$$\frac{\partial u(x,t)}{\partial t} = \frac{\partial^2 u(x,t)}{\partial x^2}, \quad 0 < x < 1, \; 0 < t \leq t_0, \tag{1}$$

$$u(x,0) = 0; \quad 0 \leq x \leq 1, \tag{2}$$

$$u(0,t) = h(t), \quad t \in [0, t_0], \tag{3}$$

$h(t) \in W_2^2[0, t_0], \; \|h(t)\|_{W_2^2}^2 = \int_0^{t_0} h(t)^2 dt + \int_0^{t_0} |h''(t)|^2 dt,$

$$h(0) = h'(0) = h(t_0) = h'(t_0) = 0, \tag{4}$$

$$\int_0^{t_0} h(t)^2 dt + \int_0^{t_0} |h''(t)|^2 dt \leq r_1^2, \tag{5}$$

$r_1$ – is a known number.

$$u(1,t) = 0, \quad t \in [0, t_0], \tag{6}$$

Consider a classical solution $U(x,t)$ of the problem (1-6), i.e. $u(x,t) \in C([0,1] \times [0, t_0]) \cap C^{2,1}((0,1) \times (0, t_0))$.

The existence and uniqueness of the solution follows from the theorem formulated in [3], p. 190. The solution of the problem (1-6) has the form:

$$u(x,t) = (1-x)h(t) + \sum_{n=1}^{\infty} v_n(t) \sin \pi n x, \tag{7}$$

$$v_n(t) = -\frac{2}{\pi n} \int_0^t e^{-(\pi n)^2(t-\tau)} h'(\tau) d\tau. \tag{8}$$

**The solution $U(x,t)$ smoothness analysis:**

It follows from (8) that

$$v_n(t) = \frac{2}{(\pi n)^3} [1 - e^{-(\pi n)^2 t}] h'(t). \tag{9}$$

It follows from (7-9) that

$$u(x,t) \in C([0,1] \times [0, t_0]). \tag{10}$$

Proceed to the examination of the function $u'_t(x,t)$ continuity. For this purpose we differentiate the common term of the series (7) relative to $t$:

$$\left[-\frac{2}{\pi n} \int_0^t e^{-(\pi n)^2(t-\tau)} h'(\tau) d\tau\right]'_t = 2h'(t) \frac{e^{-(\pi n)^2 t}}{\pi n} [1 - e^{-2(\pi n)^2 t}]. \tag{11}$$

By Abel sign and from (11) it follows that for any sufficiently small $\varepsilon > 0$ the series of the derivatives converges uniformly on the rectangle $[\varepsilon, 1-\varepsilon] \times [\varepsilon, t_0]$.

Thereby

$$u(x,t) \in C^{2,1}((0,1) \times (0, t_0]). \tag{12}$$

Using (10) and (12) we get that the solution of the problem (1-6), which has the form (7) is classical.

Using (4),(7), and(9) we get for any $x$

$$u(x, t_0) = 0. \tag{13}$$

**The statement of the inverse boundary-value problem**

Suppose, that in direct problem definition (1-6) the function $h(t)$, which defines boundary condition (3), isn't known and it should be evaluated. Therefore the additional condition is entered:
$$u(x_0,t) = f(t), x_0(0,1), t \in [0,t_0]. \qquad (14)$$
Using (7) and (14) we get
$$f(t) = (1-x_0)h(t) + \sum_{n=1}^{\infty} v_n(t)\sin \pi n x_0. \qquad (15)$$
Suppose that if $f(t) = f_0(t)$, then there is a decision $h_0(t) \in W_2^2[0,t_0]$, which satisfies conditions (4) and (5), but $f_0(t)$ it is not known and instead of it we have $f_\delta(t) \in L_2[0,t_0]$ and $\delta > 0$, for which
$$\int_0^{t_0} |f_\delta(t) - f_0(t)|^2 dt \le \delta^2. \qquad (16)$$
We need to find the approximate decision $h_\delta(t)$ and get the estimation $\|h_\delta(t) - h_0(t)\|_{L_2}$ using $f_\delta(t)$ and $\delta$.

Let's enter the linear operator $A$, which transforms the space $L_2[0,t_0]$ to $L_2[0,t_0]$ and is given by following
$$Ah(t) = -2\int_0^t K(t,\tau)h(\tau)d\tau, \text{где} \qquad (17)$$
$$K(t,\tau) = \sum_{n=1}^{\infty} \pi n e^{-(\pi n)^2(t-\tau)} \sin \pi n x_0. \qquad (18)$$
Note, that if $h_0(t) \in W_2^2[0,t_0]$ and condition (4) is realized, the inverse boundary-value problem (1-2), (5), (6), (14), (16) is equivalent to the next integral equation:
$$Ah(t) = f(t); h(t), f(t) \in L_2[0,t_0]. \qquad (19)$$
It is known that the task of the Volterr equation of the 1$^{st}$ sort solution in space $L_2[0,t_0]$ is incorrect, so for its decision we use the A.N.Tikhonov's regularization method [2].

**The 2$^{nd}$ order of the Tikhonov's regularization method**

This method consists in a reduction of the equation (19) to the variation problem, that depends on parameter $\alpha > 0$.
$$\inf\left\{\|Ah(t) - f_\delta(t)\|^2 + \alpha \int_0^{t_0} |h(t)|^2 dt + \alpha \int_0^{t_0} |h''(t)|^2 dt: h(t) \in W_2^2[0,t_0], h(0) = h(t_0) = 0\right\} (20)$$
The task (20) is equivalent to the integro-differential equation
$$A^*Ah(t) + \alpha h^{(IV)}(t) + \alpha h(t) = A^* f_\delta(t) \qquad (21)$$
$A^*$- operator, which is associated with $A$, $h(t) \in W_2^4[0,t_0]$ and $h(0) = h''(0) = h(t_0) = h''(t_0) = 0$.

It is known (see [2]) that for any $\alpha > 0$ and $f_\delta(t) \in L_2[0,t_0]$ there is an unambiguous solution $h_\delta^\alpha(t)$ of the equation (21).

The value of the regularization parameter $\alpha = \alpha(f_\delta, \delta)$ can be defined from the residual principle [4], which satisfies the equation
$$\|Ah_\delta^\alpha(t) - f_\delta(t)\|_{L_2}^2 = \delta^2. \qquad (22)$$
It is known that if $\|f_\delta(t)\|^2 > \delta^2$, then the equation (22) has the unambiguous solution $\alpha(f_\delta, \delta)$.

Thus, the approximate solution $h_\delta(t)$ of the equation (19) csn be determined by the formula:
$$h_\delta(t) = h_\delta^{\alpha(f_\delta,\delta)}(t). \qquad (23)$$

**The estimation of the error $\|h_\delta(t) - h_0(t)\|_{L_2}$**

For the error estimation the continuity module $\omega(\delta, r_1)$ will be entered:
$$\omega(\delta, r_1) = \sup\left\{\|h(t)\|_{L_2}: h(t) \in W_2^2[0,t_0], h(0) = h'(0) = h(t_0) = h'(t_0) = 0, \int_0^{t_0} h(t)^2 dt + \int_0^{t_0} |h''(t)|^2 dt \le r_1^2, \|Ah(t)\|_{L_2}^2 \le \delta^2\right\}$$
. (24)

The estimation proof is presented in work [5]
$$\|h_\delta(t) - h_0(t)\|_{L_2} \le 2\,\omega(\delta, r_1), \qquad (25)$$
where $h_\delta(t)$ is determined by (23).

Consider the expansion of the inverse problem (1), (2), (5), (6), (14) on a half-line $[t_0, \infty)$. For this purpose let's enter the functions $\bar{u}(x,t)$ and $\bar{f}(t)$ for which
$$\bar{u}(x,t) = \begin{cases} u(x,t); 0 \le x \le 1, t \in [0,t_0] \\ 0; 0 \le x \le 1, t > t_0 \end{cases} \text{и} \qquad (26)$$
$$\bar{f}(t) = \begin{cases} f(t); t \in [0,t_0] \\ 0; t > t_0 \end{cases} \text{и} \qquad (27)$$
The continuity of functions $\bar{u}(x,t)$ and $\bar{f}(t)$ follows from (10) and (13). Using (26) we get that the function $\bar{u}(x,t)$ is the solution of the problem

$$\frac{\partial \bar{u}(x,t)}{\partial t} = \frac{\partial^2 \bar{u}(x,t)}{\partial x^2}, 0 < x < 1, 0 < t, \tag{28}$$
$$\bar{u}(x,0) = 0; 0 \le x \le 1, \tag{29}$$
$$\bar{u}(x_0,t) = \bar{f}(t), t \ge 0, \tag{30}$$
$$\bar{u}(1,t) = 0; t \ge 0, \tag{31}$$

Function $\bar{h}(t)$ needs to be defined while
$$\bar{u}(0,t) = \bar{h}(t) \tag{32}$$

Let's designate as $\bar{H}$ the linear variety $L_2[0,\infty)$ such that $\bar{h}(t) \in \bar{H}$ if and only if
$$\bar{h}(t) = \begin{cases} h(t); 0 \le t \le t_0 \\ 0; t > t_0 \end{cases}, \tag{33}$$

where $h(t)$ satisfies condition (4).

Let's designate as $\bar{A}$ linear operator, which acts from $L_2[0,\infty)$ to $L_2[0,\infty)$ and which is defined on the set $\bar{H}$ by following:
$$\bar{A}\bar{h}(t) = \bar{f}(t), \tag{34}$$

where $\bar{f}(t) = \bar{u}(x_0,t)$, and $\bar{u}(x,t)$ is the solution of the problem (28), (29), (31) and (32).

For the operator $\bar{A}$ we enter the continuity module $\bar{\omega}(\delta,r_1)$:
$$\bar{\omega}(\delta,r_1) = sup\left\{\|\bar{h}(t)\|_{L_2}: \bar{h}(t) \in \bar{H}, \int_0^\infty |\bar{h}(t)|^2 dt + \int_0^\infty |\bar{h}''(t)|^2 dt \le r_1^2, \|\bar{A}\bar{h}(t)\|_{L_2}^2 \le \delta^2\right\} \tag{35}$$

Using (24), (33-35), (13) we get
$$\bar{\omega}(\delta,r_1) = \omega(\delta,r_1). \tag{36}$$

For an upper estimation of the functions $\bar{\omega}(\delta,r_1)$ we solve the problem (28-31) uing Fourier transformation on t on a half-line $[0,\infty)$.

Denote this transformation through $F_t$.

Thus the task (28-31) can be reduced to the following
$$\frac{\partial^2 \hat{u}(x,\tau)}{\partial x^2} = i\tau\hat{u}(x,\tau); 0 < x < 1, 0 < \tau, \tag{37}$$

where $\hat{u}(x,\tau) = F_t[\bar{u}(x,t)]$,
$$\hat{u}(1,\tau) = 0; \tau \ge 0 \tag{38}$$
$$\hat{u}(x_0,\tau) = \hat{f}(\tau); \tau \ge 0, \tag{39}$$

where $\hat{f}(\tau) = F_t[\bar{f}(t)]$.

The solution of the problem (37) has the form
$$\hat{u}(x,\tau) = B(\tau)e^{\mu_0 x\sqrt{\tau}} + C(\tau)e^{-\mu_0 x\sqrt{\tau}}, \tau \ge 0, \tag{40}$$

where $\mu_0 = \frac{1}{\sqrt{2}}(1+i)$. $B(\tau)$ and $C(\tau)$ need to be defined.

It follows from (38) that
$$B(\tau)e^{\mu_0\sqrt{\tau}} + C(\tau)e^{-\mu_0\sqrt{\tau}} = 0, \tau \ge 0. \tag{41}$$

It follows from (39) that
$$B(\tau)e^{\mu_0 x_0\sqrt{\tau}} 0 + C(\tau)e^{-\mu_0 x_0\sqrt{\tau}} = \hat{f}(\tau), \tau \ge 0. \tag{42}$$

Using (41) and (42) we get
$$B(\tau) = -\frac{e^{-\mu_0\sqrt{\tau}}}{2 sh\mu_0(1-x_0)\sqrt{\tau}}\hat{f}(\tau); C(\tau) = \frac{e^{\mu_0\sqrt{\tau}}}{2 sh\mu_0(1-x_0)\sqrt{\tau}}\hat{f}(\tau); \tag{43}$$

It follows from (40-43) that
$$\hat{\bar{A}}\hat{h}(\tau) = \frac{sh\mu_0(1-x_0)\sqrt{\tau}}{sh\mu_0\sqrt{\tau}}\hat{h}(\tau) = \hat{f}(\tau), \tag{44}$$

where $\hat{h}(\tau) \in F_t[\bar{H}]$, a $\hat{f}(\tau) \in L_2[0,\infty)$

It follows from conditions (4) and (5) that
$$\int_0^\infty \sqrt{1+\tau^4}\hat{h}_0(\tau)d\tau \le r_1^2. \tag{45}$$

The operator $\hat{A}$ can be dilatated on all space $L_2[0,\infty)$ without designation changing:
$$\hat{A}\hat{h}(\tau) = \frac{sh\mu_0(1-x_0)\sqrt{\tau}}{sh\mu_0\sqrt{\tau}}\hat{h}(\tau) = \hat{f}(\tau), \tag{46}$$

where $\hat{h}(\tau)$ и $\hat{f}(\tau) \in L_2[0,\infty)$

It follows from (46) that $\hat{A}$ is the injective linear bounded operator.

The relation (45) defines the the addition operator $D$, which transforms the space $L_2[0,\infty)$ to $L_2[0,\infty)$
$$D\hat{g}(\tau) = \frac{\hat{g}(\tau)}{\sqrt{1+\tau^4}}, \tau \ge 0, \text{ и} \tag{47}$$

$$\hat{h}(\tau) = D\hat{g}(\tau). \tag{48}$$

Enter the correctness class $\ddot{M}_{r_1}$

$$\hat{M}_{r_1} = D\hat{S}_{r_1}, \tag{49}$$

where $\hat{S}_{r_1} = \hat{S}(0, r_1)$ is the sphere in the space $L_2[0, \infty)$ with the center in zero with the radius $r_1$.

Let's designate as $\bar{M}_r$ the subset of $\bar{H}$ such that $\bar{h}(t) \in M_{r_1} \Leftrightarrow$

$$\int_0^{t_0} |\bar{h}(t)|^2 dt + \int_0^{t_0} |\bar{h}''(t)|^2 dt \leq r_1^2, \tag{50}$$

Using (47-49) and (50) we get

$$\ddot{M}_{r_1} \supset F_t[\bar{M}_{r_1}] \tag{51}$$

Enter the continuity module $\hat{\omega}(\delta, r_1)$ of the operator $\hat{A}$ on the spase $\ddot{M}_{r_1}$.

$$\hat{\omega}(\delta, r_1) = \left\{ \|\hat{h}(\tau)\|_{L_2} : \hat{h}(\tau) \in \hat{M}_{r_1}, \|\hat{A}\hat{h}(\tau)\| \leq \delta \right\}. \tag{52}$$

Using (50-52), (35) and using the izometry of the transformation $F_t$ we can get

$$\omega(\delta, r_1) \leq \hat{\omega}(\delta, r_1). \tag{53}$$

The estimation, that follows from (25), (36) and (53) has the form

$$\|h_\delta(t) - h_0(t)\|_{L_2} \leq 2\hat{\omega}(\delta, r_1). \tag{54}$$

Let's proceed to the estimation of the functions $\hat{\omega}(\delta, r_1)$. For this purpose we estimate the function $\frac{|sh\mu_0\sqrt{\tau}|}{|sh\mu_0(1-x_0)\sqrt{\tau}|}$. As this function is limited on any segment we get that there is the number $r_2$ such that

$$sup_{\tau \in [0,2]} \frac{|sh\mu_0\sqrt{\tau}|}{|sh\mu_0(1-x_0)\sqrt{\tau}|} \leq r_2, \tag{55}$$

when $\tau \geq 2$

$$\frac{|sh\mu_0\sqrt{\tau}|}{|sh\mu_0(1-x_0)\sqrt{\tau}|} \leq 8e^{x_0\sqrt{\frac{\tau}{2}}}, \tag{56}$$

Let's define the number $\tau_0 \geq 2$ such a way that for any $\tau \geq \tau_0$

$$e^{x_0\sqrt{\frac{\tau}{2}}} \geq r_2. \tag{57}$$

Using (57) we get for any $\tau \geq \tau_0$

$$\frac{|sh\mu_0\sqrt{\tau}|}{|sh\mu_0(1-x_0)\sqrt{\tau}|} \leq 9e^{x_0\sqrt{\frac{\tau}{2}}}, \tag{58}$$

So, when we have $\tau \geq \tau_0$

$$\frac{r_1}{\sqrt{2\tau^2}} \leq \frac{r_1}{\sqrt{1+\tau^4}}, \tag{59}$$

If $\tau_0^2 \leq e^{x_0\sqrt{\frac{\tau_0}{2}}}$, then from (58) and (59) it follows that when $\bar{\tau}$ has the form

$$\bar{\tau} = \frac{1}{2x_0^2} \ln^2\left(\frac{r_1}{9\delta}\right), \tag{60}$$

using (60) and using the theorem, proved in [6] when $\bar{\tau} \geq \tau_0$

$$\hat{\omega}(\delta, r_1) \leq \frac{r_1}{\sqrt{1+\bar{\tau}^4}}, \tag{61}$$

Or $\hat{\omega}(\delta, r_1) \sim \left[\ln\left(\frac{r_1}{9\delta}\right)\right]^{-4}$.

So, it follows from (60), (61) and (54) that for rather small values $\delta$ we get

$$\|h_\delta(t) - h_0(t)\|_{L_2} \leq \frac{r_1}{\sqrt{1 + \frac{1}{16x_0^2}[\ln\left(\frac{r_1}{9\delta}\right)]^8}}$$